\input amstex


\def\b1{\text{\bf 1}}

\def\BZ{{\Bbb Z}}

\def\CD{{\Cal D}}

\def\CT{{\Cal T}}

\def\fD{{\frak D}}


\parskip=6pt

\documentstyle{amsppt}
\document
\NoBlackBoxes


\centerline{\bf Vertex algebras and elliptic genus}

\bigskip

\centerline{Arbeitstagung, Bonn, June 1999}

\bigskip

\centerline{Vadim Schechtman}
\bigskip

\bigskip

\bigskip

This is a report on a joint work with V. Gorbounov, F. Malikov and 
A. Vaintrob, cf. [MSV], [MS1], [MS2], [GMS]. 

{\bf 1.} Let $X$ be a smooth complex variety (algebraic or analytic). 
One can associate with $X$ a canonically defined gerbe (champ of groupoids) 
$\fD^{ch}_X$ called {\it the gerbe of chiral differential operators} on $X$. 
For an open $U\subset X$, the objects of the groupoid $\fD_X(U)$ are sheaves 
of vertex algebras over $X$, called {\it sheaves of chiral differential operators} on $X$. The gerbe $\fD_X$ is bounded by the lien 
isomorphic to the sheaf $\Omega^{2, cl}_X$ of closed two-forms. 

The equivalence class
$$
c(\fD_X)\in H^2(X,\Omega^{2, cl}_X)
$$
is equal to 
$$
2ch_2(\CT_X):=c_1^2(\CT_X)-2c_2(\CT_X)
\eqno{(1.1)}
$$
where $c_i(\CT_X)\in H^i(X,\Omega^{i, cl}_X)$ are the Chern classes. 

For example, (1.1) vanishes for products of curves, or for flag spaces 
$G/B$.  

{\bf 2.} Let $E$ be a vector bundle over $X$. Generalizing the previous construction, one can define the gerbe $\fD_{\Lambda E}$ of chiral 
differential operators on the exterior algebra $\Lambda E$, also 
bounded by the sheaf $\Omega^{2, cl}_X$.  

Its equivalence class is equal to 
$$
c(\fD_{\Lambda E})=2ch_2(\CT_X)-2ch_2(E)
\eqno{(2.1)}
$$

For example, (2.1) vanishes for $E=\Omega^1_X$ or $\CT_X$. 

{\bf 3.} Consider the case $E=\Omega^1_X$. Let us denote the 
gerbe $\fD_{\Lambda E}$ of chiral differential operators 
on the de Rham algebra of differential forms by $\fD_\Omega$. 
The groupoid $\fD_\Omega(X)$ is non-empty. 

In fact, one can define certain 
{\it canonical} object $\CD_\Omega$ of $\fD_\Omega(X)$. 
It was introduced in [MSV], and called {\it chiral de Rham complex} there. 
Recently A. Beilinson gave a nice coordinate-free construction of 
$\CD_\Omega$ and some of its generalizations in the language 
of chiral algebras (in the sense of Beilinson - Drinfeld). 

$\CD_\Omega$ is a sheaf of conformal vertex superalgebras, $\BZ$-graded 
by {\it fermionic charge}. For a compact $X$, the cohomology algebra  
$$
H^*(X,\CD_\Omega)=\oplus_{p=0}^n\oplus_{q=-\infty}^\infty 
\oplus_{i=0}^\infty\ H^p(X,\CD_{\Omega, i}^q)
\eqno{(3.1)}
$$
is called the {\it chiral Hodge cohomology} of $X$ and has some nice 
properties. Here $\CD_{\Omega, i}^q$ denotes the component of 
conformal weight $i$ and fermionic charge $q$, $n=\dim(X)$.  

The space (3.1) is a conformal vertex superalgebra. It is a $N=2$ 
superalgebra if $X$ is Calabi-Yau (i.e. when the canonical bundle is trivial). 
We have the "Poincar\'e duality": 
$$
H^p(X,\CD^q_{\Omega, i})^*=H^{n-p}(X,\CD_{\Omega, i}^{n-q})
\eqno{(3.2)}
$$
As was noted in [BL], if $X$ is Calabi-Yau then the generating series 
$$
\chi_X(q, y):= y^{n/2}\sum_{a, b, i}\ (-1)^{a+b} \dim H^a(X,\CD_{\Omega, i}^b) 
y^b q^i 
\eqno{(3.3)}
$$
coincides with the Ochanine --- Witten elliptic genus of $X$.

\bigskip 

\centerline{\bf References}

\bigskip
\bigskip

[BL] L.~Borisov, A.~Libgober, Elliptic genera and applications to Mirror 
symmetry, math/9906126.   

[GMS] V.~Gorbounov, F.~Malikov, V.~Schechtman, Gerbes of chiral differential 
operators, math.AG/9906117. 

[MSV] F.~Malikov, V.~Schechtman, A.~Vaintrob, Chiral de Rham complex,
{\it Comm. Math. Phys.} (1999), to appear; math.AG/9803041.  

[MS1] F.~Malikov, V.~Schechtman, Chiral de Rham complex. II, 
{\it D.B.~Fuchs' 60-th Anniversary volume} (1999), to appear; 
math.AG/9901065. 

[MS2] F.~Malikov, V.~Schechtman, Chiral Poincar\'e duality, math.AG/9905008. 


\bigskip

\bigskip

Department of Mathematics, University of Glasgow, 
15 University Gardens, Glasgow G12 8QW, UK;\ 
vs\@maths.gla.ac.uk,\ vadik\@mpim-bonn.mpg.de

\enddocument